\documentclass[12pt]{article}
\usepackage{amsmath,amssymb,amsthm}

\title{Ramanujan's Harmonic Number Expansion into Negative Powers of a Triangular Number}

\author{Mark B. Villarino\\
Depto.\ de Matem\'atica, Universidad de Costa Rica,\\
2060 San Jos\'e, Costa Rica}

\date{July 20, 2007}

\newtheorem{lemma}{Lemma}
\newtheorem{thm}{Theorem}

\newtheorem{prop}{Proposition}
\numberwithin{equation}{section}

\makeatletter
\def\section{\@startsection{section}{1}{\z@}{-3.5ex plus -1ex minus
			  -.2ex}{2.3ex plus .2ex}{\large\bf}}
\def\subsection{\@startsection{subsection}{2}{\z@}{-3.25ex plus -1ex
			  minus -.2ex}{1.5ex plus .2ex}{\normalsize\bf}}
\renewcommand{\@dotsep}{200} 
\makeatother

\newcommand{\dis}{\displaystyle}
\newcommand{\fatt}{\rule[-18pt]{0pt}{0pt}\rule[19pt]{0pt}{0pt}}

\renewcommand{\geq}{\geqslant}  
\renewcommand{\leq}{\leqslant}  

\newcommand{\half}{{\mathchoice{\thalf}{\thalf}{\shalf}{\shalf}}}
\newcommand{\shalf}{{\scriptstyle\frac{1}{2}}} 
\newcommand{\thalf}{\tfrac{1}{2}} 

\newcommand{\hideqed}{\renewcommand{\qed}{}} 

\newcommand{\nth}{\ensuremath{n^{\mathrm{th}}}} 

\newcommand{\word}[1]{\quad\mbox{#1}\quad} 

\begin{document}

\maketitle

\begin{abstract}
An algebraic transformation of the \textsc{DeTemple--Wang}
half-integer approximation to the harmonic series produces the general
formula and error estimate for the \textsc{Ramanujan} expansion for
the $n$th harmonic number into negative powers of the $n$th
triangular number. We also discuss the history of the
\textsc{Ramanujan} expansion for the $n$th harmonic number as well as
sharp estimates of its accuracy, with complete proofs, and we compare
it with other approximative formulas.
\end{abstract}

\tableofcontents


\section{Introduction}


\subsection{The Harmonic Series}

In $1350$, \textsc{Nicholas Oresme} proved that the celebrated
\emph{\textbf{Harmonic Series}},
\begin{equation}
\label{HS}
\boxed{1 + \frac{1}{2} + \frac{1}{3} +\cdots+ \frac{1}{n} +\cdots,}
\end{equation}
is \emph{divergent}. He actually proved a more precise result. If the
\nth\ partial sum of the harmonic series, today called the \nth\
\emph{\textbf{harmonic number}}, is denoted by symbol $H_n$:
\begin{equation}
\label{Hn}
\boxed{H_n := 1 + \frac{1}{2} + \frac{1}{3} +\cdots+ \frac{1}{n},}
\end{equation}
then Oresme proved that the inequality
\begin{equation}
\label{Or }
H_{2^k} > \frac{k+1}{2}
\end{equation}
holds for $k = 2,3,\dots$. This inequality gives an estimate of the
\emph{speed} of divergence.

Almost four hundred years passed until \textsc{Leonhard Euler}, in
1755 \cite{Bromwich} applied the Euler--Maclaurin sum formula to find
the famous standard Euler \emph{\textbf{asymptotic}} expansion
for~$H_n$,
\begin{align}
\label{E}
\begin{split}
H_n := \sum_{k=1}^n \frac{1}{k}
&\sim \ln n + \gamma + \frac{1}{2n} - \frac{1}{12n^2}
+ \frac{1}{120n^4} - [\cdots]
\\
&= \ln n + \gamma - \sum_{k=1}^\infty \frac{B_k}{n^k} 
\end{split}
\end{align}
where $B_k$ denotes the $k^{\mathrm{th}}$ Bernoulli number and
$\gamma := 0.57721\cdots$ is Euler's constant. This gives a complete
answer to the speed of divergence of $H_n$ in powers of $\frac{1}{n}$.

Since then many mathematicians have contributed other approximative
formulas for $H_n$ and have studied their velocity of divergence. We
will present a detailed study of such a formula stated by
\textsc{Ramanujan}, with complete proofs, as well as of some related
formulas.


\subsection{Ramanujan's Formula}

Entry~9 of Chapter~38 of B.~\textsc{Berndt}'s edition of
Ramanujan's Notebooks \cite[p.~521]{Berndt} reads,
\begin{quote}
``\textit{Let $m := \dfrac{n(n + 1)}{2}$, where $n$ is a positive
integer. Then, as $n$ approaches infinity,}
\begin{align}
\sum_{k=1}^n \frac{1}{k}
&\sim \frac{1}{2} \ln(2m) + \gamma + \frac{1}{12m} - \frac{1}{120m^2}
+ \frac{1}{630m^3} - \frac{1}{1680m^4} + \frac{1}{2310m^5} 
\nonumber \\
&\quad  - \frac{191}{360360m^6} + \frac{29}{30030m^7}
- \frac{2833}{1166880m^8} + \frac{140051}{17459442m^9} - [\cdots]
\text{.''}
\label{Ram}
\end{align}
\end{quote}

We note that $m := \frac{n(n + 1)}{2}$ is the $n$th
\emph{\textbf{triangular number}}, so that Ramanujan's expansion of
$H_n$ is into\emph{ powers of the reciprocal of the \nth\ triangular
number.}

Berndt's proof simply verifies (as he himself explicitly notes) that
Ramanujan's expansion coincides with the standard Euler
expansion~\eqref{E}.

However, Berndt does not give the \emph{general formula} for the
coefficient of $\frac{1}{m^k}$ in Ramanujan's expansion, nor does he
prove that it is an \textit{asymptotic} series in the sense that the
error in the value obtained by stopping at any particular stage in
Ramanujan's series is less than the next term in the series. Indeed we
have been unable to find \emph{any} error analysis of Ramanujan's
series.

We will prove the following theorem.

\begin{thm}
\label{rthm}
For any integer $p \geq 1$ define
\begin{equation}
\label{R_p}
\fbox{$\dis
R_p := \frac{(-1)^{p-1}}{2p\cdot 8^p}
\biggl\{ 1 + \sum_{k=1}^p \binom{p}{k} (-4)^{k} B_{2k}(\half) \biggr\}
$}
\end{equation}
where $B_{2k}(x)$ is the Bernoulli polynomial of order~$2k$. Put
\begin{equation}
\label{m}
\fbox{$\dis
m := \frac{n(n + 1)}{2}
$}
\end{equation}
where $n$ is a positive integer. Then, for every integer $r \geq 1$,
there exists a $\Theta_r$, $0 < \Theta_r < 1$, for which the following
equation is true:
\begin{equation}
\label{expand}
\fbox{$\dis
1 +  \frac{1}{2} + \frac{1}{3} +\cdots+ \frac{1}{n}
= \frac{1}{2} \ln(2m) + \gamma + \sum_{p=1}^r \frac{R_p}{m^p}
+ \Theta_r \cdot \frac{R_{r+1}}{m^{r+1}} \,.
$}
\end{equation}
\end{thm}

We observe that the formula for $R_p$ can be written symbolically
as follows: 
\begin{equation}
\label{symR}
R_p = - \frac{1}{2p} \biggl( \frac{4B^2 - 1}{8} \biggr)^p,
\end{equation}
where we write $B_{2m}(\half)$ in place of $B^{2m}$ after carrying out
the above expansion.

We will also trace the history of Ramanujan's expansion as well and
discuss the relative accuracy of his approximation when compared to
other approximative formulas proposed by mathematicians.


\subsection{History of Ramanujan's Formula}

In 1885, two years before Ramanujan was born, \textsc{Ces\`aro}
\cite{Cesaro} proved the following.

\begin{thm}
For every positive integer $n\geq 1$ there exists a number $c_n$,
$0 < c_n < 1$, such that the following approximation is valid:
$$
\boxed{H_n = \dis \frac{1}{2}\ln(2m) + \gamma + \frac{c_n}{12m} \,.}
\eqno \qed
$$
\hideqed
\end{thm}

This gives the first two terms of Ramanujan's expansion, with an error
term. The method of proof, different from ours, does not lend itself
to generalization. We believe Ces\`aro's paper to be the first
appearance in the literature of Ramanujan's expansion.

Then, in 1904, \textsc{Lodge}, in a very interesting paper
\cite{Lodge}, which later mathematicians inexplicably (in our opinion)
ignored, proved a version of the following two results.

\begin{thm}
\label{L1}
For every positive integer $n$, define the quantity $\lambda_n$ by
the following equation:
\begin{equation}
\label{Lodge1}
\fbox{$\dis
1 + \frac{1}{2} + \frac{1}{3} +\cdots+ \frac{1}{n}
:= \frac{1}{2} \ln(2m) + \gamma + \frac{1}{12m + \frac{6}{5}}
+ \lambda_n
\,$.}
\end{equation}
Then
$$
0 < \lambda_n < \frac{19}{25200m^3}.
$$  
In fact,
$$
\lambda_n = \frac{19}{25200m^3} - \rho_n,
\word{where}  0 < \rho_n < \dfrac{43}{84000m^4}.
$$ 
The constants $\dfrac{19}{25200}$ and $\dfrac{43}{84000}$ are the best
possible.
\qed
\end{thm}

\begin{thm}
\label{L2}
For every positive integer $n$, define the quantity $\Lambda_n$ by
the following equation:
\begin{equation}
\label{Lodge2}
\boxed{1 + \frac{1}{2} + \frac{1}{3} +\cdots+ \frac{1}{n}
=: \frac{1}{2} \ln(2m) + \gamma + \frac{1}{12m + \Lambda_n} \,.}
\end{equation}
Then
$$
\Lambda_n = \frac{6}{5} - \frac{19}{175m} + \frac{13}{250m^2}
 - \frac{\delta_n}{m^3},
$$
where $0 < \delta_n < \dfrac{187969}{4042500}$. The constants in the
expansion of $\Lambda_n$ all are the best possible.
\qed
\end{thm}

These two theorems appeared, in much less precise form and \emph{with
no error estimates}, in Lodge~\cite{Lodge}. Lodge gives some numerical
examples of the error in the approximative equation
$$
H_n \approx \frac{1}{2} \ln(2m) + \gamma + \frac{1}{12m + \frac{6}{5}}
$$
in Theorem~\ref{L1}; he also presents the first two terms of
$\Lambda_n$ from Theorem~\ref{L2}. An asymptotic error estimate for
Theorem~\ref{L1} (with the incorrect constant $\frac{1}{150}$ instead
of $\frac{1}{165\frac{15}{19}}$) appears as Exercise~19 on page~460 in
\textsc{Bromwich}~\cite{Bromwich}.

\vspace{6pt}

\emph{Theorem~\ref{L1} and Theorem~\ref{L2} are immediate corollaries
of Theorem~\ref{rthm}}.

\vspace{6pt}

The next appearance of the expansion of $H_n$, into powers of the
reciprocal of the \nth\ triangular number,
$m = \dfrac{1}{\frac{n(n + 1)}{2}}$, is Ramanujan's own
expansion~\eqref{Ram}.


\subsection{Sharp Error Estimates}

Mathematicians have continued to offer alternate approximative
formulas to Euler's. We cite the following formulas, which appear in
order of increasing accuracy.

\begin{center}
\vspace{6pt}
\begin{tabular}{c|l|c|c}
\hline
No.& Approximative Formula for $H_n$& Type & Asymptotic Error Estimate
\\ \hline
$1$
& $\ln n + \gamma + \dfrac{1}{2n}$
& overestimates
& $\dfrac{1}{12n^2}$\fatt
\\ \hline
$2$
& $\ln n + \gamma + \dfrac{1}{2n + \frac{1}{3}}$
& underestimates
& $\dfrac{1}{72n^3}$\fatt 
\\ \hline
$3$
& $\ln\sqrt{n(n + 1)} + \gamma + \dfrac{1}{6n(n+1) + \frac{6}{5}}$
& overestimates
& $\dfrac{1}{165\frac{15}{19} [n(n + 1)]^3}$\fatt  
\\ \hline
$4$
& $\ln(n + \half) + \gamma + \dfrac{1}{24(n + \half)^2 + \frac{21}{5}}$
& overestimates
& $\dfrac{1}{389 \frac{781}{2071} (n + \half)^6}$\fatt
\\
\hline\hline
\end{tabular}
\vspace{6pt}
\end{center}

Formula~1 is the original Euler approximation, and it
\emph{overestimates} the true value of $H_n$ by terms of order
$\dfrac{1}{12n^2}$.
	
Formula~2 is the \textsc{T\'oth--Mare} approximation, see~\cite{TM},
and it \emph{underestimates} the true value of $H_n$ by terms of order
$\dfrac{1}{72n^3}$.

Formula~3 is the \textsc{Ramanujan--Lodge} approximation, and it
\emph{overestimates} the true value of $H_n$ by terms of order
$\dfrac{19}{3150[n(n + 1)]^3}$, see~\cite{Vill}.

Formula~4 is the \textsc{DeTemple--Wang} approximation, and it
\emph{overestimates} the true value of $H_n$ by terms of order
$\dfrac{2071}{806400 (n + \half)^6}$, see~\cite{D}.

In 2003, \textsc{Chao-Ping Chen} and \textsc{Feng Qi} \cite{CQ} gave a
proof of the following sharp form of the T\'oth--Mare approximation.

\begin{thm}
For any natural number $n \geq 1$, the following inequality is valid:
\begin{equation}
\fbox{$\dis
\frac{1}{2n + \frac{1}{1 - \gamma} - 2}
\leq H_n - \ln n - \gamma  < \frac{1}{2n + \frac{1}{3}}
\,$.}
\end{equation}
The constants $\dfrac{1}{1 - \gamma} - 2 = .3652721\cdots$ and
$d\frac{1}{3}$ are the best possible, and equality holds only for
$n = 1$.
\qed
\end{thm}

The first \emph{statement} of this theorem had been announced ten
years earlier by the editors of the ``Problems'' section of the
\emph{American Mathematical Monthly}, \textbf{99} (1992), p.~685, as
part of a commentary on the solution of Problem E~3432, but they did
not publish the proof. So, the first published proof is apparently
that of Chen and~Qi.

In this paper we will prove \emph{new and sharp forms} of the
Ramanujan--Lodge approximation and the DeTemple--Wang approximation.

\begin{thm}[Ramanujan--Lodge]
\label{Ram-Lo}
For any natural number $n \geq 1$, the following inequality is valid:
\begin{equation}
\fbox{$\dis
\frac{1}{6n(n + 1) + \frac{6}{5}} < H_n - \ln\sqrt{n(n + 1)} - \gamma
\leq \frac{1}{6n(n + 1)
+ \frac{1}{1 - \gamma - \ln\sqrt{2}}-12}
\,$.}
\end{equation}
The constants
$\dfrac{1\ln 2}{1 - \gamma - \ln\sqrt{2}}-12
= 1.12150934\cdots$ and $\dfrac{6}{5}$ are the best possible, and
equality holds only for $n = 1$.
\end{thm}

\begin{thm}[DeTemple--Wang]
\label{DeT-Wang}
For any natural number $n \geq 1$, the following inequality is valid:
\begin{equation}
\fbox{$\dis
\frac{1}{24(n + \half)^2 + \frac{21}{5}}
\leq H_n - \ln(n + \half) - \gamma
< \frac{1}{24(n + \half)^2
+ \frac{1}{1 - \ln\frac{3}{2} - \gamma}-54}
\,$.}
\end{equation}
The constants
$\dfrac{1}{1 - \ln\frac{3}{2} - \gamma}-54
= 3.73929752\cdots$ and $\dfrac{21}{5}$ are the best possible, and
equality holds only for $n = 1$.
\end{thm}

DeTemple and Wang never stated this approximation to $H_n$ explicitly.
They gave the asymptotic expansion of $H_n$, cited below in
Proposition~\ref{dw}, and we developed the corresponding approximative
formulas given above.

All three theorems are corollaries of the following stronger theorem.

\begin{thm}
\label{unlabeled}
For any natural number $n \geq 1$, define $f_n$, $\lambda_n$, and
$d_n$ by
\begin{align}
H_n &=: \ln n + \gamma + \frac{1}{2n + f_n}    
\nonumber \\
\label{lambdan}
&=: \ln\sqrt{n(n + 1)} + \gamma + \frac{1}{6n(n + 1) + \lambda_n}
\\
\label{dn}
&=: \ln(n + \half) + \gamma + \frac{1}{24(n + \half)^2 + d_n},
\end{align}
respectively. Then for any natural number $n \geq 1$ the sequence
$\{f_n\}$ is \textbf{monotonically decreasing} while the sequences
$\{\lambda_n\}$ and $\{d_n\}$ are \textbf{monotonically increasing}.
\end{thm}

Chen and Qi \cite{CQ} proved that the sequence $\{f_n\}$
\emph{decreases} monotonically. In this paper we will use their
techniques to prove the monotonicity of the sequences $\{\lambda_n\}$
and~$\{d_n\}$.


\section{Proof of the sharp error estimates}


\subsection{A few Lemmas}

Our proof is based on inequalities satisfied by the \textbf{digamma}
function $\Psi(x)$,
\begin{equation}
\fbox{$\dis
\Psi(x) := \frac{d}{dx} \ln\Gamma(x)
\equiv \frac{\Gamma'(x)}{\Gamma(x)}
\equiv - \gamma - \frac{1}{x} + x \sum_{n=1}^\infty \frac{1}{n(x + n)}
\,$,}
\end{equation}
which is the generalization of $H_n$ to the real variable $x$ since
$\Psi(x)$ and $H_n$ satisfy the equation \cite[(6.3.2), p.~258]{AS}:
\begin{equation}
\label{chiH}
\Psi(n + 1) = H_n - \gamma.
\end{equation}

\begin{lemma}
\label{lemma-one}
For every $x > 0$ there exist numbers $\theta_x$ and $\Theta_x$, with
$0 < \theta_x < 1$ and $0 < \Theta_x < 1$, for which the following
equations are true:
\begin{align}
\label{chi}
\Psi(x + 1)
&= \ln x + \frac{1}{2x} - \frac{1}{12x^2} + \frac{1}{120x^4}
- \frac{1}{252x^6} + \frac{1}{240x^8} \theta_x,
\\[2\jot]
\label{chi'}
\Psi'(x + 1)
&= \frac{1}{x} - \frac{1}{2x^2} + \frac{1}{6x^3} - \frac{1}{30x^5}
+ \frac{1}{42x^7} - \frac{1}{30x^9} \Theta_x.
\end{align}
\end{lemma}

\begin{proof}
Both formulas are well known. See, for example,
\cite[pp.~124--125]{Ed}.
\end{proof}

\begin{lemma}
\label{lemma-two}
The following inequalities are true for $x > 0$:
\begin{align}
\frac{1}{3x^2} - \frac{1}{3x^3} + \frac{4}{15x^4} - \frac{1}{5x^5}
+ \frac{10}{63x^6} - \frac{1}{7x^7}
&< 2\Psi(x + 1) - \ln\{x(x + 1)\}
\nonumber \\
&< \frac{1}{3x^2} - \frac{1}{3x^3} + \frac{4}{15x^4} - \frac{1}{5x^5}
+ \frac{10}{63x^6},
\label{lema2.1}
\\[2\jot]
\frac{2}{3x^3} - \frac{1}{4x^4} + \frac{16}{15x^5} - \frac{1}{x^6}
+ \frac{20}{21x^7} - \frac{1}{x^8}
&< \frac{1}{x} + \frac{1}{x + 1} - 2\Psi'(x + 1)
\nonumber \\
&< \frac{2}{3x^3} - \frac{1}{4x^4} + \frac{16}{15x^5} - \frac{1}{x^6}
+ \frac{20}{21x^7}.
\label{lema2.2}
\end{align}
\end{lemma}

\begin{proof}
The inequalities \eqref{lema2.1} are an immediate consequence of
\eqref{chi} and the Taylor expansion of
$$
- \ln x(x + 1)
= - 2\ln x - \ln\biggl( 1 + \frac{1}{x} \biggr)
= 2\ln\biggl( \frac{1}{x} \biggr) - \frac{1}{x} + \frac{1}{2x^2}
- \frac{1}{3x^3} + [\cdots]
$$
which is an alternating series with the property that its sum is
bracketed by two consecutive partial sums.

For \eqref{lema2.2} we start with \eqref{chi'}. We conclude that
$$
\frac{1}{2x^2} - \frac{1}{6x^3} + \frac{1}{30x^5} - \frac{1}{36x^7}
< \frac{1}{x} - \Psi'(x + 1)
< \frac{1}{2x^2} - \frac{1}{6x^3} + \frac{1}{30x^5}.
$$
Now we multiply to all three components of the inequality by~$2$ and
add $\dfrac{1}{x + 1} - \dfrac{1}{x}$ to them.
\end{proof}

\begin{lemma}
\label{lemma-three}
The following inequalities are true for $x > 0$:
\begin{align*}
\frac{1}{(x + \half)} - \frac{1}{x}
&+ \frac{1}{2x^2} - \frac{1}{6x^3} + \frac{1}{30x^5} - \frac{1}{42x^7}
< \frac{1}{x + \frac{1}{2}} - \Psi'(x + 1)
\\
&< \frac{1}{(x + \half)} - \frac{1}{x} + \frac{1}{2x^2}
- \frac{1}{6x^3} + \frac{1}{30x^5} \,,
\\
\frac{1}{24x^2} - \frac{1}{24x^3}
&+ \frac{23}{960x^4} - \frac{1}{160x^5} - \frac{11}{8064x^6}
- \frac{1}{896x^7}
< \Psi(x + 1) - \ln(x + \half)
\\
&< \frac{1}{24x^2} - \frac{1}{24x^3} + \frac{23}{960x^4}
- \frac{1}{160x^5} - \frac{11}{8064x^6} - \frac{1}{896x^7}
+ \frac{143}{30720x^8} \,.
\end{align*}
\end{lemma}

\begin{proof}
Similar to the proof of Lemma~\ref{lemma-two}.
\end{proof}


\subsection{Proof for the Ramanujan--Lodge approximation}

\begin{proof}[Proof of Theorem~\ref{unlabeled} for $\{\lambda_n\}$]
We solve \eqref{lambdan} for $\lambda_n$ and use \eqref{chiH} to
obtain
$$
\lambda_n = \frac{1}{\Psi(n + 1) - \ln\sqrt{n(n + 1)}} - 6n(n + 1).
$$
Define
$$
\Lambda_x := \frac{1}{2\Psi(x + 1) - \ln x(x + 1)} - 3x(x + 1),
$$
for all $x > 0$. Observe that $2\Lambda_n = \lambda_n$.

\vspace{6pt}

\emph{We will show that that the derivative $\Lambda_x' > 0$ for
$x > 28$.}  Computing the derivative we obtain
$$
\Lambda_x' = \frac{\frac{1}{x} + \frac{1}{x + 1} - \Psi'(x + 1)}
{\{ 2\Psi(x + 1) - \ln x(x + 1)\}^2} - (6x + 3),
$$
and therefore 
$$
\{2\Psi(x + 1) - \ln x(x + 1)\}^2 \Lambda_x'
= \frac{1}{x} + \frac{1}{x + 1} - \Psi'(x + 1)
- (6x + 3) \{2\Psi(x + 1) - \ln x(x + 1)\}^2.
$$
By Lemma~\ref{lemma-two}, this is greater than
\begin{align*}
\frac{2}{3x^3} &- \frac{1}{4x^4} + \frac{16}{15x^5} - \frac{1}{x^6}
+ \frac{20}{21x^7} - \frac{1}{x^8}
- (6x + 3) \biggl\{ \frac{1}{3x^2} - \frac{1}{3x^3} + \frac{4}{15x^4}
- \frac{1}{5x^5} + \frac{10}{63x^6} \biggr\}^2
\\[\jot]
&= \frac{798 x^5 - 21693 x^4 - 3654 x^3 + 231 x^2 + 1300 x - 2500}
{33075 x^{12}}
\\[\jot]
&= \frac{(x - 28)(798 x^4 + 651 x^3 + 14574 x^2 + 408303 x
 + 11433784) + 320143452}{33075 x^{12}}
\end{align*}
(by the remainder theorem), which is obviously \emph{positive} for
$x > 28$.  Thus, the sequence $\{\Lambda_n\}$, $n \geq 29$, is
strictly \emph{increasing}. Therefore, so is the 
sequence~$\{\lambda_n\}$.

For $n = 1,2,3,\dots,28$, we compute $\lambda_n$ directly:
\begin{align*}
\lambda_1 &=  1.1215093 &
\lambda_2 &=  1.1683646 &
\lambda_3 &=  1.1831718 &
\lambda_4 &=  1.1896217 
\\
\lambda_5 &=  1.1929804 & 
\lambda_6 &=  1.1949431 &
\lambda_7 &=  1.1961868 & 
\lambda_8 &=  1.1970233 
\\
\lambda_9 &=  1.1976125 &
\lambda_{10} &=  1.1980429 &
\lambda_{11} &=  1.1983668 &
\lambda_{12} &=  1.1986165 
\\
\lambda_{13} &=  1.1988131 &
\lambda_{14} &=  1.1989707 &
\lambda_{15} &=  1.1990988 &
\lambda_{16} &=  1.1992045 
\\
\lambda_{17} &=  1.1992926 &
\lambda_{18} &=  1.1993668 &
\lambda_{19} &=  1.1994300 &
\lambda_{20} &=  1.1994842 
\\
\lambda_{21} &=  1.1995310 &
\lambda_{22} &=  1.1995717 &
\lambda_{23} &=  1.1996073 &
\lambda_{24} &=  1.1996387 
\\
\lambda_{25} &=  1.1996664 &
\lambda_{26} &=  1.1996911 &
\lambda_{27} &=  1.1997131 &
\lambda_{28} &=  1.1997329 .
\end{align*}
Therefore, the sequence $\{\lambda_n\}$, $n \geq 1$, is a
\emph{strictly increasing sequence}.

Moreover, in Theorem~\ref{L1}, we proved that
$$
\lambda_n = \frac{6}{5} - \Delta_n,
$$
where
$0 < \Delta_n < \dfrac{38}{175n(n + 1)}$. Therefore
$$
\lim_{n\rightarrow\infty}\lambda_n = \frac{6}{5}.
\eqno \qed
$$
\hideqed
\end{proof}


\subsection{Proof for the DeTemple--Wang Approximation}

\begin{proof}[Proof of Theorem~\ref{unlabeled} for $\{d_n\}$]
Following the idea in the proof of the Lodge--Ramanujan approximation,
we solve \eqref{dn} for $d_n$ and define the corresponding
real-variable version. Let
$$
d_x := \frac{1}{\Psi(x + 1) - \ln(x + \half)} - 24(x + \half)^2.
$$
We compute the derivative, ask when is it \emph{positive}, clear the
denominator and observe that we have to solve the inequality:
$$
\biggl\{ \frac{1}{x + \half} - \Psi'(x + 1) \biggr\}
- 48(x + \half) \bigl\{ \Psi(x + 1) - \ln(x + \half) \bigr\}^2 > 0.
$$
By Lemma~\ref{lemma-three}, the left hand side of this inequality is
\begin{align*}
&> \frac{1}{x + \half} - \frac{1}{x} + \frac{1}{2x^2} - \frac{1}{6x^3}
+ \frac{1}{30x^5} - \frac{1}{42x^7}
\\
&\qquad  - 48(x + \half)
\biggl( \frac{1}{24 x^2} - \frac{1}{24 x^3} + \frac{23}{960 x^4}
- \frac{1}{160 x^5} - \frac{11}{8064 x^6} - \frac{1}{896 x^7}
+ \frac{143}{30720 x^8} \biggr)^2
\end{align*}
for all $x > 0$. This last quantity is equal to
$$
\frac{
\begin{aligned}
&( - 9018009 - 31747716 x - 14007876 x^2 + 59313792 x^3
+ 11454272 x^4 - 129239296 x^5
\\
&\quad + 119566592 x^6 + 65630208 x^7 - 701008896 x^8 - 534417408 x^9
+ 178139136 x^{10})
\end{aligned}
}{17340825600 x^{16} (1 + 2 x)} \,.
$$
The denominator is evidently positive for $x > 0$ and the numerator
can be written in the form
$$
p(x)(x - 4) + r,
$$
where 
\begin{align*}
p(x) &= 548963242092 + 137248747452 x + 34315688832 x^2
 + 8564093760 x^3 + 2138159872 x^4
\\
&\qquad + 566849792 x^5 + 111820800 x^6 + 11547648 x^7 + 178139136 x^8
+ 178139136 x^9,
\end{align*}
with remainder $r = 2195843950359$.

Therefore, the numerator is clearly positive for $x > 4$, and
therefore, the derivative $d_x^{\,\prime}$ is also positive for
$x > 4$. Finally,
\begin{align*}
d_1 &= 3.73929752\cdots   \\
d_2 &= 4.08925414\cdots   \\  
d_3 &= 4.13081174\cdots   \\
d_4 &= 4.15288035\cdots
\end{align*}
Therefore $\{d_n\}$ is an \emph{increasing} sequence for $n \geq 1$.

Now, if we expand the formula for $d_n$ into an asymptotic series in
powers of $\dfrac{1}{n + \half}$, we obtain
$$
d_n \sim \frac{21}{5} - \frac{1400}{2071(n + \half)} + \cdots
$$
(this is an immediate consequence of Proposition~\ref{dw} below) and
we conclude that
$$
\lim_{n\to\infty} d_n = \frac{21}{5}.
\eqno \qed
$$
\hideqed
\end{proof}


\section{Proof of the general Ramanujan--Lodge expansion}

\begin{proof}[Proof of Theorem~\ref{Ram-Lo}]
Our proof is founded on the half-integer approximation to $H_n$ due to
DeTemple and Wang~\cite{D}:

\begin{prop}
\label{dw}
For any positive integer $r$ there exists a $\theta_r$, with
$0 < \theta_r < 1$, for which the following equation is true:
\begin{equation}
\label{DW}
H_n = \ln (n + \half) + \gamma
+ \sum_{p=1}^r \frac{D_p}{(n + \half)^{2p}}
+ \theta_r \cdot \frac{D_{r+1}}{(n + \half)^{2r+2}},
\end{equation}
where 
\begin{equation}
\label{dp}
D_p := - \frac{B_{2p}(\half)}{2p},
\end{equation}
and where $B_{2p}(x)$ is the Bernoulli polynomial of order~$2p$.
\end{prop}

Since $(n + \half)^2 = 2m + \frac{1}{4}$, we obtain
\begin{align*}
\sum_{p=1}^r \frac{D_p}{(n + \half)^{2p}}
&= \sum_{p=1}^r \frac{D_p}{(2m)^p \bigl(1 + \frac{1}{8m}\bigr)^p}
= \sum_{p=1}^r \frac{D_p}{(2m)^p} \biggl(1 + \frac{1}{8m}\biggr)^{-p}
\\
&= \sum_{p=1}^r \frac{D_p}{(2m)^p} \sum_{k=0}^\infty \binom{-p}{k}
\frac{1}{8^k m^k}
\\
&= \sum_{p=1}^r \frac{D_p}{2^p} \sum_{k=0}^\infty (-1)^k 
\binom{k + p - 1}{k} \frac{1}{8^k} \cdot \frac{1}{m^{p+k}}
\\
&= \sum_{p=1}^r \biggl\{ \sum_{s=0}^{p-1} \frac{D_s}{2^s}
(-1)^{p-s} \binom{p - 1}{p - s} \frac{1}{8^{p-s}} \biggr\}
\cdot \frac{1}{m^p} + E_r.
\end{align*}

Substituting the right hand side of the last equation into the right
hand side of \eqref{DW} we obtain
\begin{equation}
\label{DW1}
H_n = \ln(n + \half) + \gamma + \sum_{p=1}^r
\biggl\{ \sum_{s=0}^{p-1} \frac{D_s}{2^s} (-1)^{p-s}
\binom{p - 1}{p - s} \frac{1}{8^{p-s}} \biggr\} \cdot \frac{1}{m^p}
+ E_r + \theta_r \cdot \frac{D_{r+1}}{(n + \half)^{2r+2}}.
\end{equation}

Moreover,
\begin{align*}
\ln (n + \half)
&= \frac{\ln(n + \half)^2}{2} = \frac{1}{2} \ln(2m + \tfrac{1}{4})
\\
&= \frac{1}{2} \ln(2m) + \frac{1}{2} \ln\biggl(1 + \frac{1}{8m}\biggr)
\\
&= \frac{1}{2} \ln(2m) + \frac{1}{2} \sum_{l=1}^\infty (-1)^{l-1}
\frac{1}{l\,8^lm^l}.
\end{align*}

Substituting the right-hand side of this last equation into
\eqref{DW1}, we obtain
\begin{align*}
H_n &= \frac{1}{2} \ln(2m) + \frac{1}{2} \sum_{l=1}^r (-1)^{l-1}
\frac{1}{l\,8^lm^l} + \gamma + \sum_{p=1}^r \biggl\{ 
\sum_{s=0}^{p-1} \frac{D_s}{2^s} (-1)^{p-s} \binom{p - 1}{p - s}
\frac{1}{8^{p-s}} \biggr\} \cdot \frac{1}{m^p}
\\
&\qquad + \epsilon_r + E_r
+ \theta_r \cdot \frac{D_{r+1}}{(n + \half)^{2r+2}}
\\
&= \frac{1}{2} \ln(2m) + \gamma + \sum_{p=1}^r \biggl\{
(-1)^{p-1} \frac{1}{2p\,8^p} + \sum_{s=0}^{p-1} \frac{D_s}{2^s}
(-1)^{p-s} \binom{p - 1}{p - s} \frac{1}{8^{p-s}} \biggr\}
\cdot \frac{1}{m^p}
\\
&\qquad + \epsilon_r + E_r
+ \theta_r \cdot \frac{D_{r+1}}{(n + \half)^{2r+2}}.
\end{align*}

Therefore, we have obtained Ramanujan's expansion into powers of
$\frac{1}{m}$, and the coefficient of $\frac{1}{m^p}$ is
\begin{equation}
\label{rp1}
R_p = (-1)^{p-1} \frac{1}{2p\,8^p} + \sum_{s=0}^{p-1}
\frac{D_s}{2^s} (-1)^{p-s} \binom{p - 1}{p - s} \frac{1}{8^{p-s}} \,.
\end{equation}

But,
\begin{align*}
\frac{D_s}{2^s} (-1)^{p-s} \binom{p - 1}{p - s} \frac{1}{8^{p-s}}
&= - \frac{B_{2s}(\half)/2s}{2^s} (-1)^{p-s} \binom{p - 1}{p - s}
\frac{1}{8^{p-s}}
\\
&= (-1)^{p-s-1} \frac{B_{2s}(\half)}{2s\,2^s} \binom{p - 1}{p - s}
\frac{1}{8^{p-s}},
\end{align*}
and therefore
\begin{align*}
R_p &= (-1)^{p-1} \frac{1}{2p\,8^p} + \sum_{s=0}^{p-1} \frac{D_s}{2^s}
(-1)^{p-s} \binom{p - 1}{p - s} \frac{1}{8^{p-s}}
\\
&= (-1)^{p-1} \frac{1}{2p\,8^p} + \sum_{s=0}^{p-1} (-1)^{p-s-1}
\frac{B_{2s}(\half)}{2s\,2^s} \binom{p - 1}{p - s} \frac{1}{8^{p-s}}
\\
&= (-1)^{p-1} \biggl\{ \frac{1}{2p\,8^p} + \sum_{s=1}^p (-1)^s
\frac{B_{2s}(\half)}{2s\,2^s} \binom{p - 1}{p - s} \frac{1}{8^{p-s}}
\biggr\}
\\
&= (-1)^{p-1} \biggl\{ \frac{1}{2p\,8^p} + \sum_{s=1}^p (-1)^s
\frac{B_{2s}(\half)}{2\cdot 2^s} \cdot \frac{1}{p} \binom{p}{s}
\frac{1}{8^{p-s}} \biggr\}
\\
&= \frac{(-1)^{p-1}}{2p\,8^p} \biggl\{ 1 + \sum_{s=1}^p
\binom{p}{s} (-4)^s B_{2s}(\half) \biggr\}.
\end{align*}

Therefore, the formula for $H_n$ takes the form
\begin{equation}
\label{Hn-bis}
H_n = \frac{1}{2} \ln(2m) + \gamma + \sum_{p=1}^r
\frac{(-1)^{p-1}}{2p\,8^p} \biggl\{ 1 + \sum_{s=1}^p \binom{p}{s}
(-4)^s B_{2s}(\half) \biggr\} \cdot \frac{1}{m^p} + \mathcal{E}_r,
\end{equation}
where
\begin{equation}
\label{error}
\mathcal{E}_r := \epsilon_r + E_r
+ \theta_r \cdot \frac{D_{r+1}}{(n + \half)^{2r+2}}.
\end{equation}
         
We see \emph{that \eqref{Hn-bis} is the Ramanujan expansion with the
general formula} as given in the statement of the theorem, while
\eqref{error} is a form of the \emph{error term}.
      
We will now estimate the error, \eqref{error}.

To do so, we will use the fact that the sum of a convergent
alternating series, whose terms (taken with positive sign) decrease
monotonically to zero, is equal to any partial sum \emph{plus a
positive fraction of the first neglected term (with sign)}.
      
Thus,
$$
\epsilon_r := \sum_{l=r+1}^\infty (-1)^{l-1} \frac{1}{2l\,8^lm^l}
= \alpha_r (-1)^r \frac{1}{2(r + 1) 8^{r+1} m^{r+1}},
$$
where $0 < \alpha_r < 1$. 
      
Moreover,
\begin{align*}
E_r &:= \frac{D_2}{2^1} \sum_{k=r}^\infty (-1)^k \binom{k}{k} 
\frac{1}{8^k} \cdot \frac{1}{m^{1+k}} + \frac{D_4}{2^2}
\sum_{k=r-1}^\infty (-1)^k \binom{k + 1}{k} \frac{1}{8^k}
\cdot \frac{1}{m^{2+k}} +\cdots
\\
&\qquad + \frac{D_{2r}}{2^r} \sum_{k=1}^\infty (-1)^k
\binom{k + r - 1}{k} \frac{1}{8^k} \cdot \frac{1}{m^{r+k}} + \theta_r
\cdot \frac{D_{2r+2}}{(2m)^{r+1}\bigl(1 + \frac{1}{8m}\bigr)^{r+1}}
\\
&= \biggl\{ \delta_1 \frac{D_2}{2^1} (-1)^r \binom{r}{r} \frac{1}{8^r}
+ \delta_2 \frac{D_4}{2^2} (-1)^{r-1} \binom{r}{r - 1}
\frac{1}{8^{r-1}} + \cdots
\\ 
&\qquad\qquad + \delta_r \frac{D_{2r}}{2^r} (-1)^1 \binom{r}{1}
\frac{1}{8^1} + \delta_{r+1} \frac{D_{2r+2}}{2^{r+1}} \biggr\}
\frac{1}{m^{r+1}}
\\
&= \Delta_r \biggl\{ \frac{D_2}{2^1} (-1)^r \binom{r}{r} \frac{1}{8^r}
+ \frac{D_4}{2^2} (-1)^{r-1} \binom{r}{r-1} \frac{1}{8^{r-1}} + \cdots
\\
&\qquad + \frac{D_{2r}}{2^r} (-1)^1 \binom{r}{1} \frac{1}{8^1}
+ \frac{D_{2r+2}}{2^{r+1}} \biggr\} \frac{1}{m^{r+1}},
\end{align*}
where $0 < \delta_k < 1$ for $k = 1,2,\dots,r+1$ and
$0 < \Delta_r < 1$. Thus, \emph{the error is equal to}
\begin{align*}
\mathcal{E}_r
&= \Theta_r \cdot \biggl\{ (-1)^r \frac{1}{2(r + 1)8^{(r+1)}}
+ \sum_{q=1}^{r+1} \frac{D_{2q}}{2^q} (-1)^{r - q + 1}
\binom{r}{r - q + 1} \frac{1}{8^{r-q+1}} \biggr\} \frac{1}{m^{r+1}}
\\
&= \Theta_r \cdot R_{r+1},
\end{align*}
by \eqref{R_p}, where $0 < \Theta_r < 1$, which is of the required
form. This completes the proof.
\end{proof}

The origin of Ramanujan's formula is mysterious.  Berndt notes that in his remarks.  Our analysis of it is \emph{a posteriori} and, although it is full and complete, it does not shed light on how Ramanujan came to think of his expansion.   It would also be interesting to develop an expansion for $n!$ into powers of $m$, a new \textsc{Stirling} expansion, as it were.



\begin{thebibliography}{11}

\bibitem{AS}
M. Abramowitz and I. A. Stegun,
\textit{Handbook of Mathematical Functions},
Dover, New York, 1965.

\bibitem{Berndt}
B. Berndt,
\textit{Ramanujan's Notebooks, Volume 5}, 
Springer, New York, 1998.

\bibitem{Bromwich}
T. J. l'A. Bromwich,
\textit{An Introduction to the Theory of Infinite Series},
Chelsea, New York, 1991.

\bibitem{Cesaro}
E. Ces\`aro,
``Sur la serie harmonique'',
\textit{Nouvelles Annales de Math\'ematiques} (3) \textbf{4} (1885),
295--296.

\bibitem{CQ}
Ch.-P. Chen and F. Qi,
``The best bounds of the \nth\ harmonic number'',
\emph{Global J. Math. and Math. Sci.} 2 (2006), accepted.
``The best lower and upper bounds of harmonic sequence'',
RGMIA Research Report Collection 6 (2003), no. 2, Article 14.
``The best bounds of harmonic sequence'',
arXiv:math.CA/0306233, Jiaozuo, Henan, China, 2003.

\bibitem{D}
D. DeTemple and S.-H. Wang,
``Half-integer approximations for the partial sums of the harmonic
series'',
\emph{J. Math. Anal. Appl.}, 160 (1991), 149--156. 

\bibitem{Ed}
J. Edwards,
\textit{A Treatise on the Integral Calculus}, vol.~2,
Chelsea, New York, 1955.

\bibitem{Lodge}
A. Lodge,
``An approximate expression for the value of
$\dis 1 + \frac{1}{2} + \frac{1}{3} +\cdots+ \frac{1}{r}$,''
\textit{Messenger of Mathematics} \textbf{30} (1904), 103--107.

\bibitem{TM}
L. T\'oth, and S. Mare,
Problem E~3432,
\emph{Amer. Math. Monthly}, \textbf{98} (1991), 264. 

\bibitem{Vill}
M. Villarino,
``Ramanujan's approximation to the \nth\ partial sum of the harmonic
series'',
arXiv:math.CA/0402354, San Jos\'e, 2004.

\bibitem{Vill2}
M. Villarino,
``Best bounds for the harmonic numbers'',
arXiv:math.CA/0510585, San Jos\'e, 2005.

\end{thebibliography}
\end{document}